\documentclass[a4paper]{article}%
\usepackage{graphicx}
\usepackage{geometry}
\usepackage{amsmath}
\usepackage{amssymb}
\usepackage{amsfonts}
\usepackage{amsthm,amscd}%
\providecommand{\U}[1]{\protect\rule{.1in}{.1in}}
\def\figurename{Figure}
\makeatletter
\renewcommand{\fnum@figure}[1]{\figurename~\thefigure.}
\makeatother
\def\tablename{Table}
\makeatletter
\renewcommand{\fnum@table}[1]{\tablename~\thetable.}
\makeatother
\def \bop {\noindent\textbf{Proof. }}
\def \eop {\hbox{}\nobreak\hfill
\vrule width 2mm height 2mm depth 0mm
\par \goodbreak \smallskip}
\newtheorem{theorem}{Theorem}[section]
\newtheorem{lemma}[theorem]{Lemma}
\newtheorem{corollary}[theorem]{Corollary}

\theoremstyle{definition}
\newtheorem{definition}[theorem]{Definition}

\theoremstyle{remark}
\newtheorem{remark}[theorem]{Remark}
\numberwithin{equation}{section}

\begin{document}

\title{ On optimal control of forward backward stochastic differential
equations\thanks{Partially supported by French-Algerian Scientific Program PHC
Tassili 13 MDU 887} }
\author{\textbf{F. Baghery}$^{\mathbf{1}}$\textbf{, \ \ N. Khelfallah}$^{\mathbf{2}},$
\ \ \textbf{B. Mezerdi }$^{\mathbf{2}},$ \ \textbf{I. Turpin}$^{\mathbf{1}}%
$\textbf{ }\\$^{1}$Universit\'{e} de Valenciennes, LAMAV-ISTV2, Mont Houy, \\59313 Valenciennes, Cedex 9, France\thanks{E-mails: Fouzia Baghery
(fouzia.baghery@univ-valenciennes.fr); Isabelle Turpin
(Isabelle.Turpin@univ-valenciennes.fr)}\\$^{2}$Universit\'{e} de Biskra, Laboratoire de Math\'{e}matiques \\Appliqu\'{e}es, B.P 145 Biskra (07000) Alg\'{e}rie\thanks{E-mails: Nabil
Khelfallah (nabilkhelfallah@yahoo.fr); Brahim Mezerdi (bmezerdi@yahoo.fr)}}
\maketitle

\begin{abstract}
We consider a control problem where the system is driven by a decoupled as
well as a coupled forward-backward stochastic differential equation.\ We prove
the existence of an optimal control in the class of relaxed controls, which
are measure-valued processes, generalizing the usual strict controls. The
proof is based on some tightness properties and weak convergence on the space
$\mathcal{D}$ of c\`{a}dl\`{a}g functions, endowed with the Jakubowsky
S-topology. Moreover, under some convexity assumptions, we show that the
relaxed optimal control is realized by a strict control.

\end{abstract}



\section{Introduction}

In this paper, we investigate the existence of optimal controls, for systems
driven by forward-backward stochastic differential equations (FBSDEs), of the
form:%
\begin{equation}
\left\{
\begin{array}
[c]{lll}%
dX_{t} & = & b\left(  t,X_{t},Y_{t},U_{t}\right)  dt+\sigma\left(
t,X_{t},Y_{t},U_{t}\right)  dW_{t},\\
-dY_{t} & = & h\left(  t,X_{t},Y_{t},U_{t}\right)  dt-Z_{t}dW_{t}-dM_{t},\\
X_{0}=x, &  & Y_{T}=\varphi\left(  X_{T}\right)  ,
\end{array}
\right.  \label{FBSDE}%
\end{equation}
where $\left(  M_{t}\right)  $ is a square integrable martingale, which is
orthogonal to the Brownian motion $\left(  W_{t}\right)  .$ The expected cost
over the time interval $\left[  0,T\right]  $ is given by%

\begin{equation}
J\left(  U_{.}\right)  =E\left[  \psi\left(  X_{T}\right)  +g\left(
Y_{0}\right)  +\int\nolimits_{0}^{T}l\left(  t,X_{t},Y_{t},U_{t}\right)
dt\right]  . \label{COST}%
\end{equation}

Backward stochastic differential equations (BSDEs) have been first introduced
by Pardoux and Peng, in the seminal paper \cite{PP}. Since then, the BSDE
theory became a powerful tool in many fields, such as mathematical finance,
optimal control, semi-linear and quasi-linear partial differential equations.
When the BSDE is associated to some forward stochastic differential equation,
the system is called a forward-backward stochastic differential equation
(FBSDE).\ The earliest version of such an equation appeared in Bismut
\cite{Bi}, in the stochastic version of Pontriagin maximum principle. See
\cite{EPQ, MY} for a complete account on the subject, and the references therein.

Control problems for systems governed by BSDEs and FBSDEs model many problems
arising in financial mathematics, especially the minimization of risk measures
(El Karoui and Barrieu \cite{BaElk}, Oksendal and Sulem \cite{OkSu}), the
recursive utility problems and the portfolio optimization problems. Therefore
it becomes quite natural to investigate this kind of problems, for themselves
as a class of interesting dynamical systems and for their connections to real
life problems. Many papers have been devoted to this subject, see e.g.
\cite{BKM, DZ, EPQ2, Pe, PW, Y} and the references therein. These papers have
been concerned by various forms of the stochastic maximum principle. Existence
of optimal relaxed controls for systems driven by BSDEs has been studied for
the first time in Buckdahn and al. \cite{BLRT}, by combining probabilistic
arguments as well as PDEs techniques. Then, Bahlali and al. \cite{BGM1}
investigate a control problem, with a general cost functional by using
probabilistic tools.\ The authors suppose that the generator is linear and
assume convexity of the cost function, as well as the action space. They
showed existence of an optimal strong control, that is an optimal control
adapted to the original filtration of the Brownian motion. In a second paper
\cite{BGM2}, they proved existence of a relaxed as well as of a strict control
for a system of controlled decoupled non linear FBSDE, where the diffusion
coefficient is not controlled and the generator does not depend on the second
variable $Z$. 

Our aim is to prove existence of optimal controls for systems driven by
FBSDEs. In the first part, we suppose that our FBSDE is decoupled, that is the
forward part of the equation does not contain the backward parts $Y$ and $Z$
and the diffusion coefficient depends explicitly on the control variable. We
use the formulation by martingale problems for the forward SDE and the
Meyer-Zheng compactness criteria, to prove\ the existence result. The proof is
inspired from a technique used by Pardoux \cite{P}. Note that, our result
improves \cite{BGM2} to the case where the diffusion coefficient is
controlled, \cite{EKNJ} to FBSDEs and \cite{BLRT} to continuous confficients
in the case where the generator does not depend upon $Z$. Moreover, note that
under our assumptions (the continuity of the coefficients $b$ and $\sigma$)
there are difficulties to apply directly HJB techniques as in \cite{BLRT}, to
obtain the necessary estimates of the solution of the HJB PDE, as well as its
gradient$.$ It should be mentioned that in \cite{BLRT}$,$ the coefficients of
the forward equation are Lipschitz continuous. Moreover, our approach based on
probabilistic techniques could be used for more general FBSDEs and BSDEs,
namely non Markov BSDEs, for which PDE (HJB) techniques do not work. Of
course, it should be mentioned that in \cite{BLRT}$,$ the authors treat the
case of a generator depending explicitely upon $Z$.

In the second part, we deal with a coupled FBSDEs, where the coefficients
depend on $X$ and $Y,$ but not on the second variable $Z$, with an
uncontrolled diffusion coefficient$.$ We use Jakubowsky's S-topology and a
suitable version of the Skorokhod theorem to prove the main result. Under some
additional convexity assumption, we show that the relaxed optimal control,
which is a measure-valued process, is in fact realized as a strict control.

\section{Formulation of the problem}

We study the existence of optimal controls for systems driven by FBSDEs of the
form (\ref{FBSDE}) where the cost functional over the time interval $\left[
0,T\right]  $ is given by (\ref{COST}).

We assume that $b,$ $\sigma,$ $l,$ $h$, $g$ and $\psi$ are given mappings,
$\left(  W_{t},t\geq0\right)  $ is a standard Brownian motion, defined on some
filtered probability space $\left(  \Omega,\mathcal{F},\mathcal{F}%
_{t},P\right)  $, satisfying the usual conditions, where $\left(
\mathcal{F}_{t}\right)  $ is not necessarily the Brownian filtration. $\left(
M_{t}\right)  $ is a square integrable martingale which is orthogonal to the
Brownian motion $\left(  W_{t}\right)  $ and $X,Y,Z$ are square integrable
adapted processes. The control variable $U_{t}$, called strict control, is a
measurable, $\mathcal{F}_{t}-$ adapted process with values in some compact
metric space $K$.

The objective of the controller is to minimize this cost functional, over the
class $\mathbb{U}$ of admissible controls, that is, adapted processes with
values in the set $K$, called the action space. A control $\widehat{U}$
satisfying $J\left(  \widehat{U}\right)  =\inf\left\{  J\left(  U\right)
,U\in\mathbb{U}\right\}  $ is called optimal.

Without additional convexity conditions, an optimal control may fail to exist
in the set $\mathbb{U}$ of strict controls even in deterministic control. It
should be noted that the set $\mathbb{U}$\ is not equipped with a compact
topology. The idea is then to introduce a new class of admissible controls, in
which the controller chooses at time t, a probability measure $q_{t}(du)$ on
the control set $K$, rather than an element $U_{t}\in\mathbb{U}$. These are
called relaxed controls. It turns out that this class of controls enjoys good
topological properties. If $q_{t}(du)=\delta_{U_{t}}(du)$ is a Dirac measure
charging $U_{t}$ for each $t$, then we get a strict control as a special case.
Thus the set of strict controls may be identified as a subset of relaxed controls.

To be convinced on the fact that strict controls may not exist even in the
simplest cases, let us consider a deterministic example.

The problem is to minimize the following cost function: $J(U)=\int
\nolimits_{0}^{T}\left(  X^{U}(t)\right)  ^{2}dt$ \ over the set $\mathbb{U}$
of measurable functions $U:[0,T]\rightarrow\{-1,1\}$, where $X^{U}(t)$ denotes
the solution of $dX^{U_{.}}(t)=U(t)dt,$ $X(0)=0.$ We have $\inf_{U_{.}%
\in\mathbb{U}}J(U_{.})=0$.

Indeed, consider the following sequence of controls:

\begin{center}
$U_{n}(t)=(-1)^{k}$ if $\frac{kT}{n} \leq t\leq\frac{(k+1)T}{n} $, $0\leq
k\leq n-1$.
\end{center}

Then clearly $|X^{U_{n}}(t)|\leq1/n$ and $|J(U_{n})|\leq T/n^{2}$ which
implies that $\inf_{u\in\mathbb{U}}J(u)=0$. There is however no control $U$
such that $J(U_{.})=0$. If this would have been the case, then for every $t$,
$X^{U_{.}}(t)=0$. This in turn would imply that $U_{t}=0$, which is
impossible. The problem is that the sequence $(U_{n})$ has no limit in the
space of strict controls. This limit, if it exists, will be the natural
candidate for optimality. If we identify $U_{n}(t)$ with the Dirac measure
$\delta_{U_{n}(t)}(du)$ and set $q^{n}(dt,du)=\delta_{U_{n}(t)}(du).dt$, we
get a measure on $[0,T]\times K$. Then $(q^{n}(dt,du))_{n}$ converges weakly
to $(T/2)dt$\textperiodcentered$\lbrack\delta_{-1}+\delta_{1}](du)$. This
suggests that the set\ of strict controls is too narrow and should be embedded
into a wider class enjoying compactness properties. The idea of a relaxed
control is to replace the $K$-valued process $(U_{t})$ with a $\mathcal{P}%
(K)$-valued process $(q_{t})$, where $\mathcal{P}(K)$ is the space of
probability measures equipped with the topology of weak convergence.

Let $\mathbb{V}$ be the set of Radon measures on $[0,T]\times K,$ whose
projections on $[0,T]$ coincide with the Lebesgue measure $dt$. Equipped with
the topology of stable convergence of measures, $\mathbb{V}$ is a compact
metric space (see \cite{JM}). Stable convergence is required for bounded
measurable functions $\phi(t,u),$ such that for each fixed $t\in\lbrack0,T]$,
$\phi(t,$\textperiodcentered$)$ is continuous. That is, a sequence $\left(
\mu^{n}\right)  $ in $\mathbb{V}$ converges in the stable topology to $\mu,$
if \ for every bounded measurable function $\phi:$ $[0,T]\times K\rightarrow
\mathbb{R}$ such that for each fixed $t\in\lbrack0,T]$, $\phi(t,$%
\textperiodcentered$)$ is continuous,

\begin{center}
$\int\nolimits_{0}^{T}\int\nolimits_{K}\phi(t,u).\mu^{n}(dt,du)$ converges to
$\int\nolimits_{0}^{T}\int\nolimits_{K}\phi(t,u).\mu(dt,du)$.
\end{center}

\begin{definition}
A measure-valued control on the filtered probability space $\left(
\Omega,\mathcal{F},\mathcal{F}_{t},P\right)  $ is a random variable $q$ with
values in $\mathbb{V}$, such that $q(\omega,dt,du)=dt.q(\omega,t,du)$ and
where $q(\omega,t,du)$ is progressively measurable with respect to
$(\mathcal{F}_{t})$ and such that for each $t$, $1_{(0,t]}.q$ is
$\mathcal{F}_{t}-$measurable. We denote by $\mathcal{R}$ the set of such
processes $q$.
\end{definition}

\begin{definition}
A strict control is a term $\alpha=(\Omega,\mathcal{F},\mathcal{F}_{t}%
,P,U_{t},W_{t},X_{t},Y_{t},Z_{t})$ such that

(1) $(\Omega,\mathcal{F},\mathcal{F}_{t},P)$ is a probability space equipped
with a filtration $(\mathcal{F}_{t})_{t\geq0}$ satisfying the usual conditions.

(2) $U_{t}$ is a $K$-valued process, progressively measurable with respect to
$(\mathcal{F}_{t})$.

(3) $W_{t}$ is a $(\mathcal{F}_{t},P)$- Brownian motion and $(W_{t}%
,X_{t},Y_{t},Z_{t},M_{t})$ satisfies FBSDE (\ref{FBSDE}), where $\left(
M_{t}\right)  $ is a square integrable martingale, orthogonal to $\left(
W_{t}\right)  $.
\end{definition}

\bigskip

The controls, as defined in the last definition, are called weak controls,
because of the possible change of the probability space and the Brownian
motion with $U_{t}.$

\begin{definition}
A relaxed control is a term $\alpha=(\Omega,\mathcal{F},\mathcal{F}%
_{t},P,q_{t},W_{t},X_{t},Y_{t},Z_{t},M_{t})$ such that

(1) $(\Omega,\mathcal{F},\mathcal{F}_{t},P)$ is a probability space equipped
with a filtration $(\mathcal{F}_{t})_{t\geq0}$ satisfying the usual conditions.

(2) $q$ is a measure-valued control on $\left(  \Omega,\mathcal{F}%
,\mathcal{F}_{t},P\right)  .$

(3) $W_{t}$ is a $(\mathcal{F}_{t},P)$- Brownian motion and $(W_{t}%
,X_{t},Y_{t},Z_{t},M_{t})$ satisfies the following FBSDE%
\begin{equation}
\left\{
\begin{array}
[c]{l}%
X_{t}=x+\int\nolimits_{0}^{t}\int\nolimits_{K}b\left(  s,X_{s},Y_{s},u\right)
q(s,du).ds+\int\nolimits_{0}^{t}\sigma\left(  s,X_{s},Y_{s},Z_{s}\right)
dW_{s},\\
Y_{t}=\varphi\left(  X_{T}\right)  +\int\nolimits_{t}^{T}\int\nolimits_{K}%
h\left(  s,X_{s},Y_{s},u\right)  q(s,du).ds-\int\nolimits_{t}^{T}Z_{s}%
dW_{s}-(M_{T}-M_{t}).
\end{array}
\right.
\end{equation}

where $\left(  M_{t}\right)  $ is a square integrable martingale, orthogonal
to $\left(  W_{t}\right)  .$
\end{definition}

Accordingly, the relaxed cost functional is defined by%

\begin{equation}
J\left(  q\right)  =E\left[  \psi(X_{T})+g\left(  Y_{0}\right)  +\int
\nolimits_{0}^{T}\int\nolimits_{K}l\left(  s,X_{s},Y_{s},,Z_{s},u\right)
q(s,du)ds\right]  .
\end{equation}

\begin{remark}
The appearance of the orthogonal martingale $\left(  M_{t}\right)  $ in
\ref{FBSDE} is due to the fact that the filtration associated to the optimal
control is usually larger than the Brownian filtration. By the Kunita-Watanabe
representation theorem, the conditional expectation with respect to this
filtration is a sum of a Brownian stochastic integral and an orthogonal martingale.
\end{remark}

\begin{remark}
As in classical control problems driven by It\^{o} SDEs, the cost functional
\ref{COST} \ may be defined as $J\left(  U\right)  =E\left[  \overline
{g}\left(  \overline{Y}_{0}\right)  \right]  $ where $\overline{Y}_{t}=\left(
Y_{t},Y_{t}^{k+1}\right)  ,$ with $Y_{t}^{k+1}$ the solution of the
one-dimensional BSDE
\[
\left\{
\begin{array}
[c]{c}%
-dY_{t}^{k+1}=l\left(  t,X_{t},Y_{t},U_{t}\right)  dt-\overline{Z}_{t}%
dW_{t}-d\overline{M}_{t},\\
Y_{T}^{k+1}=\psi(X_{T}),
\end{array}
\right.
\]

and $\overline{g}\left(  \overline{Y}_{0}\right)  =g\left(  Y_{0}\right)  +$
$Y_{0}^{k+1}.$

This property is known in optimal control as the equivalence between the Bolza
and Mayer problems.
\end{remark}

\textbf{Notations}

In the sequel we denote by:

\begin{itemize}
\item $\mathcal{C}([0,T];\mathbb{R}^{d})$: the space of continuous functions
from $[0,T]$ into $\mathbb{R}^{d}$, equipped with the topology of uniform convergence,

\item $\mathcal{D}([0,T];\mathbb{R}^{m})$: the Skorohod space of
c\`{a}dl\`{a}g functions from $[0,T]$ into $\mathbb{R}^{m}$, that is functions
which are continuous from the right with left hand limits endowed with the
Meyer-Zheng topology of convergence in $dt$-measure,

\item $\mathcal{S}^{2}([0,T];\mathbb{R}^{n})=\{X:[0,T]\times\Omega
\rightarrow\mathbb{R}^{n};X$ is progressively measurable and $E(\sup_{0\leq
t\leq T}|X_{t}|%
%TCIMACRO{\U{b2}}%
%BeginExpansion
{{}^2}%
%EndExpansion
)<+\infty\}$, \newline

\item $\mathcal{M}^{2}([0,T];\mathbb{R}^{k\times m})=\{Z:[0,T]\times
\Omega\rightarrow\mathbb{R}^{k\times m}$; $Z$ is progressively measurable and
$E\int_{0}^{T}|Z_{t}|^{2}dt<+\infty\}$, \newline

\item $L_{\mathcal{F}}^{2}([0,T];\mathbb{R}^{k})=\{f\left(  t,\omega\right)
:[0,T]\times\Omega\rightarrow\mathbb{R}^{k};\,\,\,\mathcal{F}_{t}$ - adapted
such that $E\left[  \int_{0}^{T}|f\left(  t,\omega\right)  |^{2}dt\right]
<+\infty\}$.
\end{itemize}

\section{ Control of decoupled FBSDEs}

In this section we deal with decoupled FBSDEs, that is the forward part does
not contain $(Y,Z)$ the solution of the backward part. We suppose that the
diffusion coefficient $\sigma$ depends explicitly on the control variable and
the driver $h$ as well as the instantaneous cost $l$ do not depend on $Z$.
More precisely, our system is governed by the following equation%

\begin{equation}
\left\{
\begin{array}
[c]{lll}%
dX_{t} & = & b\left(  t,X_{t},U_{t}\right)  dt+\sigma\left(  t,X_{t}%
,U_{t}\right)  dW_{t},\\
X\left(  0\right)  & = & x,\\
-dY_{t} & = & h\left(  t,X_{t},Y_{t},U_{t}\right)  dt-Z_{t}dW_{t}-dM_{t},\\
Y_{T} & = & \varphi\left(  X_{T}\right)  ,\label{3.1}%
\end{array}
\right.
\end{equation}
defined on some filtered probability space $\left(  \Omega,\mathcal{F}%
,\mathcal{F}_{t},P\right)  $, $\left(  W_{t}\right)  $ is an $m$-dimensional
Brownian motion and $\left(  M_{t}\right)  $ is a square integrable martingale
which is orthogonal to $\left(  W_{t}\right)  .$ The coefficients of our FBSDE
are defined as follows
\begin{align*}
b  &  :\left[  0,T\right]  \times\mathbb{R}^{d}\times K\rightarrow
\mathbb{R}^{d},\\
\sigma &  :\left[  0,T\right]  \times\mathbb{R}^{d}\times K\rightarrow
\mathbb{R}^{d\times m},\\
h  &  :\left[  0,T\right]  \times\mathbb{R}^{d}\times\mathbb{R}^{k}\times
K\rightarrow\mathbb{R}^{k}\\
\varphi &  :\mathbb{R}^{d}\rightarrow\mathbb{R}^{k}.
\end{align*}

Let us define the cost functional over $\left[  0,T\right]  $ by%

\begin{equation}
J\left(  U\right)  =E\left[  \psi(X_{T}) + g\left(  Y_{0}\right)
+\int\nolimits_{0}^{T}l\left(  t,X_{t},Y_{t},U_{t}\right)  dt\right]  ,
\end{equation}

where \
\begin{align*}
l  &  :\left[  0,T\right]  \times\mathbb{R}^{d}\times\mathbb{R}^{k}\times
K\rightarrow\mathbb{R},\\
\psi &  :\mathbb{R}^{d}\rightarrow\mathbb{R},\\
g  &  :\mathbb{R}^{k}\rightarrow\mathbb{R}.
\end{align*}

$\mathbf{(H}_{\mathbf{1}}\mathbf{)}$ Assume that the functions $b,$ $\sigma,$
$h,$ $\varphi$ are continuous and bounded.

$\mathbf{(H}_{\mathbf{2}}\mathbf{)}$ $h$ is Lipschitz in the variable $y$
uniformly in $(t,x,u),$ i.e: there exists a constant $C>0$ such that for every
$t\in\lbrack0,T]$, $u\in K$, $y\ ,y^{\prime}\in\mathbb{R}^{d},$
\[
\left\vert h\left(  t,x,y,u\right)  -h\left(  t,x,y^{\prime},u\right)
\right\vert \leq C\left\vert y-y^{\prime}\right\vert .
\]

$\mathbf{(H}_{\mathbf{3}}\mathbf{)}$ Assume that $l$, $\psi$ and $g$ are
continuous and bounded functions.

The infinitesimal generator $L$, associated with the forward part of our
equation, acting on functions $f$ in $C_{b}^{2}(R^{d};R)$, is defined by

\begin{center}
$Lf(t,x,u)=\left(  \dfrac{1}{2}\underset{i,j}{%
%TCIMACRO{\dsum }%
%BeginExpansion
{\displaystyle\sum}
%EndExpansion
}a_{ij}\dfrac{\partial^{2}f}{\partial x_{i}x_{j}}+\underset{i}{%
%TCIMACRO{\dsum }%
%BeginExpansion
{\displaystyle\sum}
%EndExpansion
}b_{i}\dfrac{\partial f}{\partial x_{i}}\right)  (t,x,u)$,
\end{center}

where $a_{ij}(t,x,u)$ denotes the generic term of the symmetric matrix
$\sigma\sigma^{\ast}(t,x,u)$.

As it is well known, weak solutions for It\^{o} SDEs are equivalent to the
existence of solutions for the corresponding martingale problem. Then one can
rewrite definition 2.1 and definition 2.2, by using the formulation of
martingale problems for the forward part. This simplifies taking limits and
does not pose the problem of the relaxation of the stochastic integral part
see \cite{BDM, BM, EKNJ}. We can define a strict control using martingale
problems as follows.

\begin{definition}
A strict control is a term $\alpha=(\Omega,\mathcal{F},\mathcal{F}_{t}%
,P,U_{t},X_{t},Y_{t},M_{t})$ such that

(1) $(\Omega,\mathcal{F},\mathcal{F}_{t},P)$ is a probability space equipped
with a filtration $(\mathcal{F}_{t})_{t\geq0}$ satisfying the usual conditions.

(2) $U_{t}$ is an $K$-valued process, progressively measurable with respect to
$(\mathcal{F}_{t})$.

(3i) $(X_{t})$ is $\mathbb{R}^{d}-$valued $\mathcal{F}_{t}-$ adapted, with
continuous paths, such that
\begin{equation}
f(X_{t})-f(x)-\int_{0}^{t}Lf(s,X_{s},U_{s}%
)ds\,\,\,\mathrm{is\,\,a\,\,P-{martingale}}. \label{FMP}%
\end{equation}

(3ii) $(Y_{t},M_{t})$ is the solution of the following backward SDE%
\begin{equation}
Y_{t}=\varphi\left(  X_{T}\right)  +\int\nolimits_{t}^{T}h\left(
s,X_{s},Y_{s},U_{s}\right)  ds-\left(  M_{T}-M_{t}\right)  , \label{BMP}%
\end{equation}
where $\left(  M_{t}\right)  $ is a square integrable $(\mathcal{F}_{t})-$
martingale$.$
\end{definition}

\begin{definition}
A relaxed control is a term $\alpha=(\Omega,\mathcal{F},\mathcal{F}%
_{t},P,q_{t},X_{t},Y_{t},M_{t})$ such that

(1) $(\Omega,\mathcal{F},\mathcal{F}_{t},P)$ is a probability space equipped
with a filtration $(\mathcal{F}_{t})_{t\geq0}$ satisfying the

usual conditions.

(2) $q$ is a measure-valued control on $\left(  \Omega,\mathcal{F}%
,\mathcal{F}_{t},P\right)  .$

(3i) $(X_{t})$ is $\mathbb{R}^{d}-$valued $\mathcal{F}_{t}-$adapted, with
continuous paths, such that
\begin{equation}
f(X_{t})-f(x)-\int_{0}^{t}\int_{K}Lf(s,X_{s}%
,u).q(s,du)ds\,\,\,\mathrm{is\,\,a\,\,P-{martingale}}.
\end{equation}

(3ii) $(Y_{t},M_{t})$ is the solution of the following backward SDE
\begin{equation}%
\begin{array}
[c]{lll}%
Y_{t} & = & \varphi\left(  X_{T}\right)  +\int\nolimits_{t}^{T}\int
\nolimits_{K}h\left(  s,X_{s},Y_{s},u\right)  .q(s,du).ds-\left(  M_{T}%
-M_{t}\right)  ,
\end{array}
\end{equation}
where $\left(  M_{t}\right)  $ is a square integrable $(\mathcal{F}_{t})-$martingale.
\end{definition}

By a slight abuse of notation, we will often denote a relaxed control by $q$
instead of specifying all the components.

The cost functional associated to a relaxed control is now defined by%
\begin{equation}
J\left(  q\right)  =E\left[  \psi(X_{T})+g\left(  Y_{0}\right)  +\int
\nolimits_{0}^{T}\int\nolimits_{K}l\left(  s,X_{s},Y_{s},u\right)
q(s,du)ds\right]  .
\end{equation}
\bigskip

The main result of this section is given by the following Theorem.

\begin{theorem}
Under assumptions $\mathbf{(H_{1})}$, $\mathbf{(H_{2})}$ and $\mathbf{(H_{3}%
)}$, the relaxed control problem has an optimal solution.
\end{theorem}

The proof is based on some auxiliary results on the tightness of the processes
under consideration and the identification of the limits.

Let $\left(  q^{n}\right)  _{n\geq0}$ be a minimizing sequence, that is
$\underset{n\rightarrow\infty}{\lim}J\left(  q^{n}\right)  =\underset
{q\in\mathcal{R}}{\inf}J\left(  q\right)  $ and let $(X^{n},Y^{n},M^{n})$ be a
solution of our FBSDE, where:

(i) $(X_{t}^{n})$ is $\mathbb{R}^{d}-$valued $\mathcal{F}_{t}-$adapted, with
continuous paths, such that

\begin{center}
$f(X_{t}^{n})-f(x)-\int\nolimits_{0}^{t}\int\nolimits_{K}Lf(s,X_{s}%
^{n},u).q^{n}(s,du).ds$ is a $P$-martingale
\end{center}

(ii) $(Y_{t}^{n},M_{t}^{n})$ is the solution of the following backward SDE

\begin{center}
$Y_{t}^{n}=\varphi\left(  X_{T}^{n}\right)  +\int\nolimits_{t}^{T}%
\int\nolimits_{K}h\left(  s,X_{s}^{n},Y_{s}^{n},u \right)  .q^{n}%
(s,du).ds-\left(  M_{T}^{n}-M_{t}^{n}\right)  ,$
\end{center}

and $\left(  M_{t}^{n}\right)  $ is a square integrable continuous
$(\mathcal{F}_{t})-$martingale.

The proof of the main result consists in proving that the sequence of
distributions of the processes $(q^{n},X^{n},Y^{n},M^{n})$ is tight for a
certain topology, on the state space and then show that, we can extract a
subsequence, which converges in law to a process $(\widehat{q},\widehat
{X},\widehat{Y},\widehat{M}),$ satisfying the same FBSDE. To complete the
proof, we show that under some regularity conditions the sequence of cost
functionals $(J(q^{n}))_{n}$ converge to $J(\widehat{q})$ which is equal to
$\underset{q\in\mathcal{R}}{\inf}J\left(  q\right)  $ and then $(\widehat
{q},\widehat{X},\widehat{Y},\widehat{M})$ is optimal.

\begin{lemma}
\ The family of relaxed controls $(q^{n})_{n}$ is tight in $\mathbb{V}$.
\end{lemma}

\bop$[0,T]\times K$ being compact, then by Prokhorov's theorem, the space
$\mathbb{V}$ of probability measures on $[0,T]\times K$ is also compact for
the topology of weak convergence. The fact that $q^{n}$, $n\geq0$ are random
variables with values in the compact set $\mathbb{V}$ yields that the family
of distributions associated to $(q^{n})_{n\geq0}$ is tight. \eop

\begin{lemma}
i) The family $(q^{n},X^{n})_{n}$ of solutions of the martingale problem is
tight on the space $\mathbb{V}\times\mathcal{C}(\left(  \left[  0,T\right]
;R^{d}\right)  $.

ii) There exists a subsequence which converges in law to $(\widehat
{q},\widehat{X}),$ whose law is a solution of the martingale problem, that is
for each $f$ $\in\mathcal{C}_{b}^{2}$, $f(\widehat{X}_{t})-f(x)-\int
\nolimits_{0}^{t}\int\nolimits_{K}Lf(s,\widehat{X}_{s},u).\widehat
{q}(s,du).ds$ is a $P$-martingale
\end{lemma}

\bop Let us give the outlines of the proof which is inspired from \cite{EKNJ},
Theorem 3.4.\newline i) Following \cite{SV}, Theorem 1.4.6, it is sufficient
to show that for each positive $f$ in $\mathcal{C}_{b}^{2},$ there exists a
constant $A_{f}$ such that: $f(X_{t})+A_{f}.t$ is a supermartingale under the
distribution $P_{n}=P^{(q^{n},X^{n})}$ on the canonical space $\mathbb{V\times
}\mathcal{C} \left(  \left[  0,T\right]  ; R^{d}\right)  $ of the couple
$(q^{n},X^{n})$. Let

\begin{center}
$A_{f}$ $=\sup$ $\left\{  |Lf(t,x,u)|; (t,x,u)\in\left[  0,T\right]
\times\mathbb{R}^{d}\times K\right\}  $.
\end{center}

$A_{f}$ is finite since the coefficients $b$ and $\sigma$ defining the
operator $\mathit{L}$ are bounded. \newline Since for each $n$, $f(X_{t}%
)-f(x)-\int\nolimits_{0}^{t}\int\nolimits_{K}Lf(s,X_{s},u).q(s,du).ds:=C_{t}%
f(x,q)$ is a $P_{n}$-martingale, then $f(X_{t})+A_{f}.t$ is a positive
supermartingale. Then $(X^{n})$ is tight in $\mathcal{C}$ endowed with the
topology of uniform convergence. \newline ii) The sequence $(q^{n},X_{n})$
being tight, then we can extract a subsequence still denoted by $(q^{n}%
,X^{n})$ which converges weakly to $(\widehat{q},\widehat{X}).$ In particular,
for every bounded $(x,q)$-continuous, $\mathcal{C}_{s}\mathbb{\otimes
}\mathbb{V}_{s}$-measurable functions ($\mathcal{C}_{s}$ and $\mathbb{V}_{s}$
are the $\sigma-$fields generated by the coordinates until for $t\leq s)$, we have

\begin{center}
$P_{n}\left[  \phi(q,x)\left(  C_{t}f(x,q)-C_{s}f(x,q)\right)  \right]  $
converges to $\widehat{P}\left[  \phi(q,x)\left(  C_{t}f(x,q)-C_{s}%
f(x,q)\right)  \right]  $
\end{center}

where $\widehat{P}$ denotes the law of the couple $(\widehat{q},\widehat{X})$
in the space $\mathbb{V\times}\mathcal{C}\left(  \left[  0,T\right]  ;
R^{d}\right)  $. \newline$C_{t}f(x,q)$ being a $P_{n}$-martingale, then
$P_{n}\left[  \phi(x,q)\left(  C_{t}f(x,q)-C_{s}f(x,q)\right)  \right]  =0$.
\newline Hence the limit $\widehat{P}\left[  \phi(q,x)\left(  C_{t}%
f(x,q)-C_{s}f(x,q)\right)  \right]  =0$ and thus the law $\widehat{P}$ of the
couple $(\widehat{q},\widehat{X})$ is a solution of the martingale problem.
\eop

\begin{lemma}
i) The sequence $(Y^{n},M^{n})$ is tight on the space $\mathcal{D}^{2}$
equipped with the Meyer-Zheng topology.

ii) There exists a subsequence still denoted by $(q^{n},X^{n},Y^{n},M^{n})$
which converges weakly to $(\widehat{q},\widehat{X},\widehat{Y},\widehat{M}),$
in the space $\mathbb{V\times}\mathcal{C}\mathbb{\times}\mathcal{D}^{2}$.
Moreover $(\widehat{q},\widehat{X},\widehat{Y},\widehat{M})$ satisfies:

For every $f$ $\in\mathcal{C}_{b}^{2}$,%

\[
f(\widehat{X}_{t})-f(x)-\int\nolimits_{0}^{t}\int\nolimits_{K}Lf(s,\widehat
{X}_{s},u).\widehat{q}(s,du).ds\text{ is a }\mathcal{F}_{t}^{\widehat
{q},\widehat{X},\widehat{Y}}-martingale
\]

\begin{equation}
\widehat{Y}_{t}=\varphi\left(  \widehat{X}_{T}\right)  +\int\nolimits_{t}%
^{T}\int\nolimits_{K}h\left(  s,\widehat{X}_{s},\widehat{Y}_{s},u\right)
\widehat{q}(s,du).ds-\left(  \widehat{M}_{T}-\widehat{M}_{t}\right)
\end{equation}
where $\left(  \widehat{M}_{t}\right)  $ is a square integrable $\mathcal{F}%
_{t}^{\widehat{q},\widehat{X},\widehat{Y}}-$martingale.
\end{lemma}

\bop i) Using standard techniques from BSDEs theory it is not difficult to
prove that \newline$E\left(  \sup_{0\leq t\leq T}\left\vert Y_{t}%
^{n}\right\vert ^{2}\right)  \leq CE\left(  \left\vert \varphi\left(
X_{T}^{n}\right)  \right\vert ^{2}+\int\nolimits_{0}^{T}\sup_{u\in
K}\left\vert h\left(  s,\widehat{X}_{s},0,u\right)  \right\vert ^{2}ds\right)
$. Then using assumption \textbf{(H}$_{\mathbf{1}}$\textbf{)} it follows that:

\begin{center}
$\underset{n}{\sup}E\left(  \sup_{0\leq t\leq T}\left\vert Y_{t}%
^{n}\right\vert ^{2}+\left\langle M^{n}\right\rangle _{T}\right)  <+\infty.$
\end{center}

Let us denote the conditional variation

\begin{center}
$V_{t}(Y^{n})=\sup E\left(
%TCIMACRO{\dsum }%
%BeginExpansion
{\displaystyle\sum}
%EndExpansion
\left\vert E \left[  \left(  Y_{t_{i+1}}^{n}-Y_{t_{i}}^{n}\right)  /F_{t_{i}%
}^{X^{n}} \right]  \right\vert \right)  , $
\end{center}

where the supremum is taken over all the partitions of the interval $\left[
0,t\right]  .$ One can easily prove that%

\[
V_{t}(Y^{n})\leq E\left(  \int\nolimits_{0}^{T}\sup_{u\in K}\left\vert
h\left(  s,X_{s}^{n},Y_{s}^{n},u\right)  \right\vert ds\right)  .
\]

The assumptions made on the coefficients ensure that

\begin{center}
$\underset{n}{\sup}\left[  V_{t}(Y^{n})+\sup_{0\leq t\leq T}E\left\vert
Y_{t}^{n}\right\vert +\sup_{0\leq t\leq T}E\left\vert M_{t}^{n}\right\vert
\right]  <+\infty.$
\end{center}

Then following \cite{P, MZ}, the sequence $(Y^{n},M^{n})$ satisfies the
Meyer-Zheng criterion for tightness of families of quasi-martingales.

ii) Since the sequence $\left(  q^{n},X^{n},Y^{n},M^{n}\right)  _{n}$ is
tight, then there exist a subsequence still denoted by $(q^{n},X^{n}%
,Y^{n},M^{n})_{n}$ which converges weakly to $(\widehat{q},\widehat
{X},\widehat{Y},\widehat{M})$ on the space $\mathbb{V}\times\mathcal{C}%
\times\mathcal{D}^{2},$ where $\mathcal{C}$ is equipped with the topology of
uniform convergence and $\mathcal{D}$ is equipped with the Meyer-Zheng
topology. By using the fact that for each $t\leq T,$ the mapping
$(q,x,y)\longrightarrow\int\nolimits_{t}^{T}\int_{K}h(t,x_{s},y_{s}%
,u).q(s,du)ds$ is continuous from $\mathbb{V}\times\mathcal{C}\times
\mathcal{D}$ into $\mathbb{R}$, one can pass to the limit in the BSDE and get

\begin{center}
$\widehat{Y}_{t}=\varphi\left(  \widehat{X}_{T}\right)  +\int\nolimits_{t}%
^{T}\int\nolimits_{K}h\left(  s,\widehat{X}_{s},\widehat{Y}_{s},u\right)
.\widehat{q}(s,du).ds-\left(  \widehat{M}_{T}-\widehat{M}_{t}\right)  .$
\end{center}

Let us show that $\widehat{M}_{t}$ and $\phi(\widehat{X}_{t})-\phi(\widehat
{X}_{s})-\int\nolimits_{0}^{t}\int\nolimits_{K}L\phi(s,\widehat{X}%
_{s},u).\widehat{q}(s,du).ds$ are both martingales with respect to the natural
filtration $\mathcal{F}_{t}=\mathcal{F}_{t}^{\widehat{q},\widehat{X}%
,\widehat{Y}}.$ For any $s,t$ such that $0\leq s\leq t\leq T$ and $\Phi_{s}$ a
bounded continuous mapping from $\mathbb{V}_{s}\mathbb{\times}\mathcal{C}%
\left(  \left[  0,s\right]  ,\mathbb{R}\right)  \times\mathcal{D}\left(
\left[  0,s\right]  ,\mathbb{R}\right)  $ and for each $\phi\in C_{b}^{2}$:

\begin{center}
$E\left[  \Phi_{s}(q^{n},X^{n},Y^{n})\left(  \phi(X_{t}^{n})-\phi(X_{s}%
^{n})-\int\nolimits_{0}^{t}\int\nolimits_{K}L\phi(s,X_{s}^{n},u).q^{n}%
(s,du).ds\right)  \right]  $ $\longrightarrow0$ as $n\longrightarrow\infty$
\end{center}

and for each $n$

\begin{center}
$E\left[  \Phi_{s}(q^{n}, X^{n},Y^{n})\left(  \int\nolimits_{0}^{\varepsilon
}\left(  M_{t+r}^{n}-M_{s+r}^{n}\right)  .dr\right)  \right]  =0,$
\end{center}

where $\mathbb{V}_{s}$ denotes the restriction of measures to the interval
$\left[  0,s\right]  .$ The fact that $(q^{n},X^{n},Y^{n})$ is weakly
convergent and $E(\underset{0\leq t\leq T}{\sup}\left\vert M_{t}%
^{n}\right\vert ^{2})$ is finite yield

\begin{center}
$E\left[  \Phi_{s}(\widehat{q},\widehat{X},\widehat{Y})\left(  \phi
(\widehat{X}_{t})-\phi(\widehat{X}_{s})-\int\nolimits_{0}^{t}\int
\nolimits_{K}L\phi(s,\widehat{X}_{s},u).\widehat{q}(s,du).ds\right)  \right]
=0$

$E\left[  \Phi_{s}(\widehat{q},\widehat{X},\widehat{Y})\left(  \int
\nolimits_{0}^{\varepsilon}\left(  \widehat{M}_{t+r}-\widehat{M}_{s+r}\right)
.dr\right)  \right]  =0.$
\end{center}

In the last equality, dividing by $\varepsilon$ and then sending it to $0$ and
by the right continuity of the martingale $\left(  \widehat{M}_{t}\right)  $,
we get

\begin{center}
$E\left[  \Phi_{s}(\widehat{q}, \widehat{X},\widehat{Y})\left(  \widehat
{M}_{t}-\widehat{M}_{s}\right)  \right]  =0.$
\end{center}

These identities are valid for all functions $\Phi_{s}$ described above and
for all $s\leq t$. Then $\widehat{M}_{t}$ and $\phi(\widehat{X}_{t}%
)-\phi(\widehat{X}_{s})-\int\nolimits_{0}^{t}\int\nolimits_{K}L\phi
(s,\widehat{X}_{s},u).\widehat{q}(s,du).ds$ are both $\mathcal{F}%
_{t}^{\widehat{q},\widehat{X},\widehat{Y}}-$martingales. \eop

\textbf{Proof of Theorem 3.5.} According to Lemma 3.6 and the assumptions
\textbf{(H}$_{\mathbf{1}}$\textbf{)-(H}$_{\mathbf{3}}$\textbf{)} we have
\begin{align*}
\underset{q\in\mathcal{R}}{\inf}J\left(  q\right)   &  =\underset
{n\rightarrow\infty}{\lim}J\left(  q^{n}\right) \\
&  =\underset{n\rightarrow\infty}{\lim}E\left[  \psi(X_{T}^{n})+g\left(
Y_{0}^{n}\right)  +\int\nolimits_{0}^{T}\int\nolimits_{K}l\left(  t,X_{t}%
^{n},Y_{t}^{n},u\right)  q_{t}^{n}\left(  du\right)  dt\right] \\
&  =E\left[  \psi(\hat{X}_{T})+g\left(  \hat{Y}_{0}\right)  +\int
\nolimits_{0}^{T}\int\nolimits_{K}l\left(  t,\hat{X}_{t},\hat{Y}_{t},u\right)
\hat{q}_{t}\left(  du\right)  dt\right]  ,
\end{align*}
which means that $\hat{q}$ is an optimal control. \eop

\begin{corollary}
Assume that $\mathbf{(H_{1})}$, $\mathbf{(H_{2})}$ and $\mathbf{(H_{3})}$
hold. Moreover assume that for every $\left(  t,x,y\right)  \in\left[
0,T\right]  \times\mathbb{R}^{d}\times\mathbb{R}^{k}$, the set
\begin{equation}
\left(  b,\sigma,h,l\right)  \left(  t,x,y,K\right) \\
:=\left\{  b_{i}\left(  t,x,u\right)  ,\left(  \sigma\sigma^{\ast}\right)
_{ij}\left(  t,x,u\right)  ,h_{j}\left(  t,x,y,u\right)  ,l\left(
t,x,y,u\right)  /u\in K,i=1,...,d,\text{ \ \ }j=1,...,k\right\}
\end{equation}
is convex in $\mathbb{R}^{d+d\times d+k+1}.$ Then, the relaxed optimal control
$\hat{q}_{t}$ has the form of a Dirac measure charging a strict control
$\hat{U}_{t}$, that is $\hat{q}_{t}\left(  du\right)  =\delta_{\hat{U}_{t}%
}\left(  du\right)  $.
\end{corollary}

\bop We put

\begin{center}
$\int_{K}h\left(  t,\hat{X}_{t},\hat{Y}_{t},u\right)  \hat{q}_{t}\left(
du\right)  :=\hat{h}\left(  t,w\right)  \in h\left(  t,x,y,\mathbb{U}\right)
,$ $\int_{K}l\left(  t,\hat{X}_{t},\hat{Y}_{t},u\right)  \hat{q}_{t}\left(
du\right)  :=\hat{l}\left(  t,w\right)  \in l\left(  t,x,y,\mathbb{U}\right)
,$ $\int_{K}b\left(  t,\hat{X}_{t}^{n},u\right)  \hat{q}_{t}\left(  du\right)
:=\widehat{b}\left(  t,w\right)  \in b\left(  t,x,\mathbb{U}\right)  ,\int
_{K}a\left(  t,\hat{X}_{t}^{n},u\right)  \hat{q}_{t}\left(  du\right)
:=\widehat{a}\left(  t,w\right)  \in a\left(  t,x,\mathbb{U}\right)  ,$ where
$a=\sigma\sigma^{\ast}.$
\end{center}

From $\mathbf{(H}_{1}\mathbf{)-(H}_{3}\mathbf{)}$ and the measurable selection
theorem (see \cite{YZ} p. 74 or \cite{EKNJ}), there is a $K$-valued,
$\mathcal{F}^{\hat{X},\hat{Y},\hat{q}}-$adapted process $\hat{U},$ such that
for every $s\in\lbrack0,\ T]$,
\begin{align*}
\left(  \hat{h},\hat{l}\right)  \left(  s,w\right)   &  =\left(  h,l\right)
\left(  s,\hat{X}\left(  s,w\right)  ,\hat{Y}\left(  s,w\right)  ,\hat
{U}\left(  s,w\right)  \right)  ,\\
\left(  \hat{b},\hat{a}\right)  \left(  s,w\right)   &  =(b,a)\left(
s,\hat{X}\left(  s,w\right)  ,\hat{U}\left(  s,w\right)  \right)  .
\end{align*}
Hence, for every $t\in\lbrack0,\ T]$ and $w\in\hat{\Omega}$, we have
\[
\int_{K}h\left(  t,\hat{X}_{t},\hat{Y}_{t},u\right)  \hat{q}_{t}\left(
du\right)  =h\left(  t,\hat{X}_{t},\hat{Y}_{t},\hat{U}_{t}\right)  ,
\]%
\[
\int_{K}l\left(  t,\hat{X}_{t},\hat{Y}_{t},u\right)  \hat{q}_{t}\left(
du\right)  =l\left(  t,\hat{X}_{t},\hat{Y}_{t},\hat{U}_{t}\right)  ,
\]
and
\[
\int_{K}b\left(  t,\hat{X}_{t},u\right)  \hat{q}_{t}\left(  du\right)
=b\left(  t,\hat{X}_{t},\hat{U}_{t}\right)  ,
\]%
\[
\int_{K}a\left(  t,\hat{X}_{t},u\right)  \hat{q}_{t}\left(  du\right)
=a\left(  t,\hat{X}_{t},\hat{U}_{t}\right)  .
\]

Then the process $(\hat{X}_{t},\hat{Y}_{t},\widehat{M}_{t})$ satisfies, for
each $t\in\lbrack0,\ T]$:

1) $(\widehat{X}_{t})$ is a $\mathbb{R}^{d}-$valued $\mathcal{F}_{t}-$adapted,
with continuous paths, such that

\begin{center}
$f(\widehat{X}_{t})-f(x)-\int\nolimits_{0}^{t} Lf(s,\widehat{X}_{s},\hat
{U}_{s})ds$ is a $P$-martingale
\end{center}

2) $(\widehat{Y}_{t},\widehat{M}_{t})$ solves the following BSDE

\begin{center}
$\widehat{Y}_{t}=\varphi\left(  \widehat{X}_{T}\right)  +\int\nolimits_{t}%
^{T}h\left(  s,\widehat{X}_{s},\widehat{Y}_{s},\hat{U}_{s}\right)  ds-\left(
\widehat{M}_{T}-\widehat{M}_{t}\right)  $
\end{center}

It follows that \ $J(\widehat{q})=J(\hat{U})$, which achieves the
proof.\textbf{ }\eop

\section{ Control of coupled FBSDEs}

In this section, we consider coupled FBSDEs, where all the coefficients depend
only on $(X,Y)$ but not on the second backward component $Z$. Moreover the
diffusion coefficient $\sigma$ does not depend on the control variable. More
precisely, the controlled FBSDE takes the form%
\begin{equation}
\left\{
\begin{array}
[c]{lll}%
dX_{t} & = & b\left(  t,X_{t},Y_{t},U_{t}\right)  dt+\sigma\left(
t,X_{t},Y_{t}\right)  dW_{t},\\
-dY_{t} & = & h\left(  t,X_{t},Y_{t},U_{t}\right)  dt-Z_{t}dW_{t}-dM_{t},\\
X_{0}=x, &  & Y_{T}=\varphi\left(  X_{T}\right)  ,
\end{array}
\right.  \label{FBSDE1}%
\end{equation}

and the cost functional is given by%

\begin{equation}
J\left(  U_{.}\right)  =E\left[  \psi(X_{T})+g\left(  Y_{0}\right)
+\int\nolimits_{0}^{T}l\left(  t,X_{t},Y_{t},U_{t}\right)  dt\right]
,\label{COST1}%
\end{equation}

where $\psi,$ $g$, $l$ satisfy (\textbf{H}$_{\mathbf{3}}$).

Assume that the coefficients \
\begin{align*}
b  &  :\left[  0,T\right]  \times\mathbb{R}^{d}\times\mathbb{R}^{k}\times
K\rightarrow\mathbb{R}^{d},\\
\sigma &  :\left[  0,T\right]  \times\mathbb{R}^{d}\times\mathbb{R}%
^{k}\rightarrow\mathbb{R}^{d\times m},\\
h  &  :\left[  0,T\right]  \times\mathbb{R}^{d}\times\mathbb{R}^{k}\times
K\rightarrow\mathbb{R}^{k},\\
\varphi &  :\mathbb{R}^{d}\rightarrow\mathbb{R}^{k},
\end{align*}
of the FBSDE \ref{FBSDE1} satisfy the following conditions:

$\mathbf{(H}_{4}\mathbf{)}$ $b,$ $\sigma,$ $h$ are bounded measurable,
Lipshitz in $(x,y)$ uniformly in $(t,u)$ and continuous in $u.$

$\mathbf{(H_{5})}$ Assume that for each admissible control $U,$ the
coefficients $b,\sigma,h,\varphi$ of the FBSDE \ref{FBSDE1} satisfy the
monotonicity conditions as in \cite{PW} Theorem 2.2, page 828.

\begin{remark}
Under assumption $\mathbf{(H_{5})}$ and for each admissible control $U,$ the
FBSDE \ref{FBSDE1} admits a unique strong solution.
\end{remark}

The main result of this section is given by the following theorem.

\begin{theorem}
The relaxed control problem defined by \ref{FBSDE1} and \ref{COST1} has an
optimal solution.
\end{theorem}

The proof is based on tightness properties of the underlying processes.

As in the last section, let $\left(  q^{n}\right)  _{n\geq0}$ be a minimizing
sequence for the relaxed control problem, that is
\[
\underset{n\rightarrow\infty}{\lim}J\left(  q^{n}\right)  =\underset
{q\in\mathcal{R}}{\inf}J\left(  q\right)  ,
\]

Let $(X^{n},Y^{n},Z^{n})$ be the unique strong solution of our FBSDE
\begin{equation}
\left\{
\begin{array}
[c]{l}%
X_{t}^{n}=x+\int\nolimits_{0}^{t}\int\nolimits_{K}b\left(  s,X_{s}^{n}%
,Y_{s}^{n},u\right)  q^{n}(s,du).ds+\int\nolimits_{0}^{t}\sigma\left(
s,X_{s}^{n},Y_{s}^{n}\right)  dW_{s}\\
Y_{t}^{n}=\varphi\left(  X_{T}^{n}\right)  +\int\nolimits_{t}^{T}%
\int\nolimits_{K}h\left(  s,X_{s}^{n},Y_{s}^{n},u\right)  q^{n}(s,du).ds-\int
\nolimits_{t}^{T}Z_{s}^{n}dW_{s},
\end{array}
\right.  \label{EDSR-RC}%
\end{equation}

defined on the natural filtration of the Brownian motion $\left(
W_{t}\right)  .$ In this case the orthogonal martingales $\left(  M_{t}%
^{n}\right)  $ disappears, due to the uniqueness of solutions.

The proof of Theorem 4.2 consists in showing that the sequence $(q^{n}%
,X^{n},Y^{n},\int\nolimits_{0}^{.}Z_{s}^{n}dW_{s})$ is tight and there exists
a subsequence converging weakly to $(\widehat{q},\widehat{X},\widehat
{Y},\widehat{N})$. Furthermore these processes satisfy the FBSDE%
\begin{equation}
\left\{
\begin{array}
[c]{l}%
\widehat{X}_{t}=x+\int\nolimits_{0}^{t}\int\nolimits_{K}b\left(  s,\widehat
{X}_{s},\widehat{Y}_{s},u\right)  \widehat{q}(s,du).ds+\int\nolimits_{0}%
^{t}\sigma\left(  s,\widehat{X}_{s},\widehat{Y}_{s}\right)  d\widehat{W}%
_{s},\\
\widehat{Y}_{t}=\varphi\left(  \widehat{X}_{T}\right)  +\int\nolimits_{t}%
^{T}\int\nolimits_{K}h\left(  s,\widehat{X}_{s},\widehat{Y}_{s},u\right)
\widehat{q}(s,du).ds-\int\nolimits_{t}^{T}\widehat{Z}_{s}d\widehat{W}%
_{s}-(\widehat{M}_{T}-\widehat{M}_{t}),
\end{array}
\right.
\end{equation}
where $\widehat{M}_{t}$ is a square integrable martingale which is orthogonal
to the Brownian motion $\widehat{W}_{t}.$

\begin{lemma}
\textit{Let $(X^{n},Y^{n},Z^{n})$ be the unique solution of equation
(\ref{EDSR-RC}). There exists a positive constant }$C$\textit{\ such that}
\begin{equation}
\underset{n}{\sup}E\left(  \underset{0\leq t\leq T}{\sup}\left\vert X_{t}%
^{n}\right\vert ^{2}+\underset{0\leq t\leq T}{\sup}\left\vert Y_{t}%
^{n}\right\vert ^{2}+\int_{0}^{T}|Z_{s}^{n}|^{2}ds\right)  \leq C .
\label{borne-Yn,Zn}%
\end{equation}

\end{lemma}

\bop Using assumption $\mathbf{(H}_{\mathbf{5}}\mathbf{)}$, it is easy to
check that
\begin{equation}
\underset{n}{\sup}\{E(\underset{0\leq t\leq T}{\sup}|X_{t}^{n}|^{2}%
)\}<\infty\label{borne-Xn}%
\end{equation}
Using the Burkholder-Davis-Gundy and Schwarz inequalities, it follows that the
local martingale $\int_{t}^{T}Y_{s}^{n}Z_{s}^{n}dW_{s}$ is uniformly
integrable. Then by It\^{o}'s formula and assumption \textbf{(H}$_{\mathbf{2}%
}$\textbf{)}, it holds that%
\[
E(|Y_{t}^{n}|^{2}+\int_{t}^{T}|Z_{s}^{n}|^{2}ds)=E\left(  |\varphi(X_{T}%
^{n})|^{2}+2\int_{t}^{T}\int_{K}\langle Y_{s}^{n},h(s,X_{s}^{n},Y_{s}%
^{n},u)\rangle q_{s}^{n}(du)ds\right)  .
\]
Hence,
\begin{align*}
E(|Y_{t}^{n}|^{2}+\int_{t}^{T}|Z_{s}^{n}|^{2}ds)\leq E(  &  |\varphi(X_{T}%
^{n})|^{2}+\int_{t}^{T}|Y_{s}^{n}|^{2}ds)\\
&  +E(\int_{t}^{T}\int_{K}|h(s,X_{s}^{n},Y_{s}^{n},u)|^{2}q_{s}^{n}(du)ds).
\end{align*}
The result follows from Gronwall's Lemma and BDG inequality

\begin{center}
$\underset{n}{\sup}E\left(  \underset{0\leq t\leq T}{\sup}|Y_{t}^{n}|^{2}%
+\int_{0}^{T}|Z_{s}^{n}|^{2}ds \right)  <\infty.$
\end{center}

\eop

\begin{lemma}
\textit{Let $(X^{n},Y^{n},Z^{n})$ be the unique solution of equation
(\ref{EDSR-RC}). Then the sequence (}$Y^{n},\int_{0}^{\cdot}Z_{s}^{n}dW_{s}%
)$\textit{\ is tight on the space }$\mathbb{D}\left(  \left[  0,T\right]
;\mathbb{R}^{k}\right)  \times\mathbb{D}\left(  \left[  0,T\right]
;\mathbb{R}^{k}\right)  $\textit{\ endowed with the Jakubowski }%
$S$\textit{-topology.}
\end{lemma}

\bop Let $0=t_{0}<t_{1}<...<t_{n}=T$, and define the conditional variation by
\[
CV\left(  Y^{n}\right)  :=\sup E\left[  \underset{i}{\sum}\left\vert E\left(
Y_{t_{i+1}}^{n}-Y_{t_{i}}^{n}\right)  \diagup\mathcal{\mathcal{F}}_{t_{i}}%
^{W}\right\vert \right]  ,
\]
where the supremum is taken over all partitions of the interval $[0,T]$. It is
proved in \cite{P} that
\[
CV\left(  Y^{n}\right)  \leq\mathbb{E}\left[  \int_{0}^{T}\int_{K}|h\left(
s,X_{s}^{n},Y_{s}^{n},u\right)  |q_{s}^{n}\left(  du\right)  ds\right]  .
\]
It follows from (\ref{borne-Yn,Zn}) that
\[
\underset{n}{\sup}\left[  {CV}\left(  Y^{n}\right)  +\underset{0\leq t\leq
T}{\sup}E\left\vert Y_{t}^{n}\right\vert +\underset{0\leq t\leq T}{\sup
}E\left\vert \int_{0}^{t}Z_{s}^{n}dW_{s}\right\vert \right]  <\infty
\]
Therefore, the sequence $\left(  Y^{n}\mathit{,}\int_{0}^{\cdot}Z_{s}%
^{n}dW_{s}\right)  $\ satisfies the Meyer-Zheng criterion \cite{MZ}. Therefore
$\mathit{(}Y^{n},\int_{0}^{\cdot}Z_{s}^{n}dW_{s})$ is tight in the Jakubowski
S-topology (see the appendix). \eop

The next lemma may be proved by standard arguments.

\begin{lemma}
Let $X_{t}^{n}$ be the forward component of equation (\ref{EDSR-RC}).Then the
sequence of processes $\left(  X^{n},W\right)  $ is tight on the space
$\mathcal{C}\left(  \left[  0,T\right]  ,\mathbb{R}^{d}\right)  \times
\mathcal{C}\left(  \left[  0,T\right]  ,\mathbb{R}^{m}\right)  ,$ endowed with
the topology of uniform convergence.
\end{lemma}

\bop\textbf{ of Theorem 4.2} From Lemmas 4.3, 4.4 and 4.5, the sequence
$\gamma^{n}=\left(  q^{n},X^{n},W,Y^{n},N^{n}\right)  $ where $N^{n}=\int
_{0}^{\cdot}Z_{s}^{n}dW_{s},$ is tight on the space $\Gamma=\mathbb{V\times
}\mathcal{C}\left(  \left[  0,T\right]  ,\mathbb{R}^{d}\right)  \times
\mathcal{C}\left(  \left[  0,T\right]  ,\mathbb{R}^{m}\right)  \mathbb{\times
}\left[  \mathbb{D}\left(  \left[  0,T\right]  ;\mathbb{R}^{k}\right)
\right]  ^{2}.$ \newline Using the Skorokhod representation theorem on the
space $\mathcal{D}$ endowed with S-topology \cite{J1} (see the appendix),
there exists a probability space $\left(  \hat{\Omega},\mathcal{\hat
{\mathcal{F}}},\mathbb{\hat{P}}\right)  $, a sequence $\hat{\gamma}%
^{n}=\left(  \hat{q}^{n},\hat{X}^{n},\hat{W}^{n},\hat{Y}^{n},\hat{N}%
^{n}\right)  $ and $\hat{\gamma}=\left(  \hat{q},\hat{X},\hat{W},\hat{Y}%
,\hat{N}\right)  $ defined on this space such that:

(i) for each $n\in\mathbb{N}$, law$\left(  \gamma^{n}\right)  =$ law$\left(
\hat{\gamma}^{n}\right)  $,

(ii) there exists a subsequence $\left(  \hat{\gamma}^{n_{k}}\right)  $\ of
$\left(  \hat{\gamma}^{n}\right)  $, still denoted $\left(  \hat{\gamma}%
^{n}\right)  $, which converges to $\hat{\gamma},\mathbb{\hat{P}}$-a.s. on the
space $\Gamma,$

(iii) the subsequence $(\hat{Y}^{n},\hat{N}^{n})$ converges to $(\hat{Y}%
,\hat{N}),$ $dt\times\mathbb{\hat{P}}-a.s$., and $(\hat{Y}_{T}^{n},\hat{N}%
_{T}^{n})$ converges to $(\hat{Y}_{T},\hat{N}_{T})$ as $n\rightarrow\infty,$
$\mathbb{\hat{P}}-a.s$.

(iv) $\underset{0\leq t\leq T}{\sup}\left\vert \hat{X}_{t}^{n}-\hat{X}%
_{t}\right\vert \rightarrow0,\ \ \ \mathbb{\hat{P}}-a.s$.

According to property (i), it follows that%
\begin{equation}
\left\{
\begin{array}
[c]{ccc}%
\hat{X}_{t}^{n} & = & x+\int_{0}^{t}\int_{K}b\left(  s,\hat{X}_{s}%
^{n},\widehat{Y}_{s}^{n},u\right)  \hat{q}_{s}^{n}\left(  du\right)
ds+\int_{0}^{t}\sigma\left(  s,\hat{X}_{s}^{n},\widehat{Y}_{s}^{n}\right)
d\hat{W}_{s}^{n}\\
\hat{Y}_{t}^{n} & = & \varphi\left(  \hat{X}_{T}^{n}\right)  +\int_{t}^{T}%
\int_{K}h\left(  s,\hat{X}_{s}^{n},\hat{Y}_{s}^{n},u\right)  \hat{q}_{s}%
^{n}\left(  du\right)  ds-\left(  \hat{N}_{T}^{n}-\hat{N}_{t}^{n}\right)  ,
\end{array}
\right.  \label{EDSR-n-chapeau}%
\end{equation}
where $\hat{N}_{t}^{n}:=\int_{t}^{T}\hat{Z}_{s}^{n}d\hat{W}_{s}^{n}$.

Combining properties (ii)-(iv), assumptions ($\mathbf{H_{2}}$)$-$%
($\mathbf{H_{5}}$) and passing to the limit in the FBSDE (\ref{EDSR-n-chapeau}%
), there exists a countable set $D\subset\left[  0,\ T\right)  $ such that
\begin{equation}
\left\{
\begin{array}
[c]{cccc}%
\hat{X}_{t} & = & x+\int_{0}^{t}\int_{K}b\left(  s,\hat{X}_{s},\widehat{Y}%
_{s},u\right)  \hat{q}_{s}\left(  du\right)  ds+\int_{0}^{t}\sigma\left(
s,\hat{X}_{s},\widehat{Y}_{s}\right)  d\hat{W}_{s}, & \ t>0,\\
\hat{Y}_{t} & = & \varphi\left(  \hat{X}_{T}\right)  +\int_{t}^{T}\int
_{K}h\left(  s,\hat{X}_{s},\hat{Y}_{s},u\right)  \hat{q}_{s}\left(  du\right)
ds-\left(  \hat{N}_{T}-\hat{N}_{t}\right)  , & \ t\in\left[  \,0,\ T\right]
\setminus D.
\end{array}
\right.  \label{EDSR-chapeau}%
\end{equation}
Since $\hat{Y}$ and $\hat{N}$ are c\`{a}dl\`{a}g, it follows that for every
$t\in\left[  0,T\right]  $,%
\[
\hat{Y}_{t}=\varphi\left(  \hat{X}_{T}\right)  +\int_{t}^{T}\int_{K}h\left(
s,\hat{X}_{s},\hat{Y}_{s},u\right)  \hat{q}_{s}\left(  du\right)  ds+\hat
{N}_{t}-\hat{N}_{T}.
\]
Since all the previous identifications of the limits can be proved similarly,
let us prove that$:$
\begin{equation}
\lim_{n\rightarrow\infty}\int_{t}^{T}\int_{K}h\left(  s,\hat{X}_{s}^{n}%
,\hat{Y}_{s}^{n},u\right)  \hat{q}_{s}^{n}\left(  du\right)  ds=\int_{t}%
^{T}\int_{K}h\left(  s,\hat{X}_{s},\hat{Y}_{s},u\right)  \hat{q}_{s}\left(
du\right)  ds. \label{limitfn}%
\end{equation}
We use properties (i), (ii), (iv), Fatou's lemma and Lemma 4.3, to show that
there exists a constant $C$ such that:%
\begin{equation}
\hat{E}(\int_{0}^{T}(|\hat{X}_{s}|^{2}+|\hat{Y}_{s}|^{2})ds)\leq C.
\label{borneXY}%
\end{equation}
On the other hand, we have
\[
\left\vert \int_{t}^{T}\int_{K}h\left(  s,\hat{X}_{s}^{n},\hat{Y}_{s}%
^{n},u\right)  \hat{q}_{s}^{n}\left(  du\right)  ds-\int_{t}^{T}\int
_{K}h\left(  s,\hat{X}_{s},\hat{Y}_{s},u\right)  \hat{q}_{s}\left(  du\right)
ds\right\vert \leq I(n)+J(n),
\]
where
\[
I(n):=\left\vert \int_{t}^{T}\int_{K}h\left(  s,\hat{X}_{s}^{n},\hat{Y}%
_{s}^{n},u\right)  \hat{q}_{s}^{n}\left(  du\right)  ds-\int_{t}^{T}\int
_{K}h\left(  s,\hat{X}_{s},\hat{Y}_{s},u\right)  \hat{q}_{s}^{n}\left(
du\right)  ds\right\vert ,
\]%
\[
J(n):=\left\vert \int_{t}^{T}\int_{K}h\left(  s,\hat{X}_{s},\hat{Y}%
_{s},u\right)  \hat{q}_{s}^{n}\left(  du\right)  ds-\int_{t}^{T}\int
_{K}h\left(  s,\hat{X}_{s},\hat{Y}_{s},u\right)  \hat{q}_{s}\left(  du\right)
ds\right\vert .
\]
Let us show that $I(n)$ converges to 0 in probability. Let $\varepsilon>0$ and
use the fact that $h$ is Lipschitz in $\left(  x,y\right)  $ to obtain,%
\begin{align*}
&  \mathbb{\hat{P}}\left\{  \left\vert \int_{t}^{T}\int_{K}h\left(  s,\hat
{X}_{s}^{n},\hat{Y}_{s}^{n},u\right)  \hat{q}_{s}^{n}\left(  du\right)
ds-\int_{t}^{T}\int_{K}h\left(  s,\hat{X}_{s},\hat{Y}_{s},u\right)  \hat
{q}_{s}^{n}\left(  du\right)  ds\right\vert >\varepsilon\right\} \\
&  \ \leq\frac{1}{\varepsilon}\hat{E}\int_{t}^{T}\int_{K}\left\vert h\left(
s,\hat{X}_{s}^{n},\hat{Y}_{s}^{n},u\right)  -h\left(  s,\hat{X}_{s},\hat
{Y}_{s},u\right)  \right\vert \hat{q}_{s}^{n}\left(  du\right)  ds\\
&  \leq\frac{C}{\varepsilon}\left[  \hat{E}\int_{t}^{T}\left\vert \hat{X}%
_{s}^{n}-\hat{X}_{s}\right\vert ds+\hat{E}\int_{t}^{T}\left\vert \hat{Y}%
_{s}^{n}-\hat{Y}_{s}\right\vert ds\right]  .
\end{align*}
Now, properties (i)-(iv) and Lemma 4.3 allow us to show that $\hat{E}\int
_{t}^{T}|\hat{X}_{s}^{n}-\hat{X}_{s}|ds+\hat{E}\int_{t}^{T}|\hat{Y}_{s}%
^{n}-\hat{Y}_{s}|ds$ tends to $0$ as $n$ tends to infinity, which yields that
$I(n)$ converges to 0 in probability.

Now let us show that $J(n)$ converges to 0 in probability. Let $R>0$ and, put
$B:=\{|\hat{X}_{s}|+|\hat{Y}_{s}|\leq R\}$ and $\bar{B}:=\Omega-B$. We have,
\[
\left\vert \int_{t}^{T}\int_{K}h\left(  s,\hat{X}_{s},\hat{Y}_{s},u\right)
\hat{q}_{s}^{n}\left(  du\right)  ds-\int_{t}^{T}\int_{K}h\left(  s,\hat
{X}_{s},\hat{Y}_{s},u\right)  \hat{q}_{s}\left(  du\right)  ds\right\vert \leq
I_{1}(n)+J_{1}(n),
\]
where
\[
I_{1}(n)=:\left\vert \int_{t}^{T}\int_{K}h\left(  s,\hat{X}_{s},\hat{Y}%
_{s},u\right)  1_{B}\,\hat{q}_{s}^{n}\left(  du\right)  ds-\int_{t}^{T}%
\int_{K}h\left(  s,\hat{X}_{s},\hat{Y}_{s},u\right)  1_{B}\,\hat{q}_{s}\left(
du\right)  ds\right\vert ,
\]%
\[
J_{1}(n):=\left\vert \int_{t}^{T}\int_{K}h\left(  s,\hat{X}_{s},\hat{Y}%
_{s},u\right)  1_{\bar{B}}\,\hat{q}_{s}^{n}\left(  du\right)  ds-\int_{t}%
^{T}\int_{K}h\left(  s,\hat{X}_{s},\hat{Y}_{s},u\right)  1_{\bar{B}}\,\hat
{q}_{s}\left(  du\right)  ds\right\vert .
\]
Since the function $(s,u)\longmapsto h\left(  s,\hat{X}_{s},\hat{Y}%
_{s},u\right)  1_{B}$ is bounded, measurable in $(s,u)$ and continuous in $u$,
we deduce by using property (ii) that $I_{1}(n)$ tends to 0 in probability as
$n$ tends to $\infty$. It remains to prove that $J_{1}(n)$ tends to 0 in
probability as $n$ tends to $\infty$. We have,
\begin{align*}
\hat{E}[J_{1}(n)]  &  =\hat{E}(|\int_{t}^{T}\int_{K}h\left(  s,\hat{X}%
_{s},\hat{Y}_{s},u\right)  1\!\!1_{\bar{B}}\,\hat{q}_{s}^{n}\left(  du\right)
ds-\int_{t}^{T}\int_{K}h\left(  s,\hat{X}_{s},\hat{Y}_{s},u\right)
1\!\!1_{\bar{B}}\,\hat{q}_{s}\left(  du\right)  ds|)\\
&  \leq\hat{E}\int_{t}^{T}\int_{K}\left\vert h\left(  s,\hat{X}_{s},\hat
{Y}_{s},u\right)  \right\vert 1\!\!1_{\bar{B}}\hat{q}_{s}^{n}\left(
du\right)  ds+\hat{E}\int_{t}^{T}\int_{K}\left\vert h\left(  s,\hat{X}%
_{s},\hat{Y}_{s},u\right)  \right\vert 1\!\!1_{\bar{B}}\hat{q}_{s}\left(
du\right)  ds\\
&  \leq\frac{C^{\prime}}{R^{2}}\hat{E}\int_{t}^{T}(\left\vert \hat{X}%
_{s}\right\vert ^{2}+\left\vert \hat{Y}_{s}\right\vert ^{2})ds.
\end{align*}
We successively pass to the limit in $n$ and $R$, to show that $\lim
_{n\rightarrow\infty}J(n)=0$ in probability.

Now, let $\mathcal{\hat{F}}_{s}:=\mathcal{F}_{s}^{\hat{X},\hat{Y},\hat{q}}$,
be the filtration generated by $(\hat{X_{r}},\,\hat{Y_{r}},\hat{q_{r}},\ r\leq
s)$ completed by $\hat{P}-$nul sets. Combining the estimates
(\ref{borne-Yn,Zn}) and Lemma 6.3 in Appendix, one can show that $\left(
\hat{N}_{t}\right)  $ is a $\mathcal{\hat{F}}_{t}$-martingale. Since $\left(
\hat{W}\right)  $ is a $\left(  \mathcal{\hat{F}}_{t},\hat{P}\right)
-$Brownian motion, then the martingale decomposition theorem yields the
existence of a process $\hat{Z}\in\mathcal{M}^{2}([t,T];\mathbb{R}^{n\times
d})$ such that
\[
\hat{N}_{t}=\int_{0}^{t}\hat{Z}_{s}d\hat{W}_{s}+\hat{M}_{t},\ \ \ \text{with}%
\ \ \ \left\langle \hat{M}_{t},\hat{W}\right\rangle _{t}=0,
\]
which implies that
\[
\hat{Y}_{t}=\varphi\left(  \hat{X}_{T}\right)  +\int_{t}^{T}\int_{K}h\left(
s,\hat{X}_{s},\hat{Y}_{s},u\right)  \hat{q}_{s}\left(  du\right)  ds-\int
_{t}^{T}\hat{Z}_{s}d\hat{W}_{s}-\left(  \hat{N}_{T}-\hat{N}_{t}\right)  .
\]
To finish the proof of Theorem 4.2, it remains to check that $\hat{q}$ is an
optimal control.

\ According to above properties (i)-(iv) and assumption \textbf{(H$_{3}$)}, we
have
\begin{align*}
\underset{q\in\mathcal{R}}{\inf}J\left(  q\right)   &  =\underset
{n\rightarrow\infty}{\lim}J\left(  q^{n}\right)  ,\\
&  =\underset{n\rightarrow\infty}{\lim}E\left[  \psi(X_{T}^{n})+g\left(
Y_{0}^{n}\right)  +\int\nolimits_{0}^{T}\int\nolimits_{K}l\left(  t,X_{t}%
^{n},Y_{t}^{n},u\right)  q_{t}^{n}\left(  du\right)  dt\right] \\
&  =\underset{n\rightarrow\infty}{\lim}\hat{E}\left[  \psi(\hat{X}_{T}%
^{n})+g\left(  \hat{Y}_{0}^{n}\right)  +\int\nolimits_{0}^{T}\int
\nolimits_{K}l\left(  t,\hat{X}_{t}^{n},\hat{Y}_{t}^{n},u\right)  \hat{q}%
_{t}^{n}\left(  du\right)  dt\right] \\
&  =\hat{E}\left[  \psi(\hat{X}_{T})+g\left(  \hat{Y}_{0}\right)
+\int\nolimits_{0}^{T}\int\nolimits_{K}l\left(  t,\hat{X}_{t},\hat{Y}%
_{t},u\right)  \hat{q}_{t}\left(  du\right)  dt\right]  .
\end{align*}
\bigskip\eop

By using the same arguments as in Corollary 3.7, one can prove the following
result on the existence of strict controls under convexity assumptions.

\begin{corollary}
Assume $\mathbf{(H_{2})}$- $\mathbf{(H_{5})}$ and that for every $\left(
t,x,y\right)  \in\left[  0,T\right]  \times\mathbb{R}^{d}\times\mathbb{R}^{k}%
$, the set
\begin{equation}
\left(  b,h,l\right)  \left(  t,x,y,K\right)  :=\left\{  b_{i}\left(
t,x,y,u\right)  ,h_{j}\left(  t,x,y,u\right)  ,l\left(  t,x,y,u\right)  /u\in
K,i=1,...,d,j=1,...,k\right\}  ,
\end{equation}
is convex and closed in $\mathbb{R}^{d+k+1}.$ Then, the relaxed optimal
control $\hat{q}_{t}$ has the form of a Dirac measure charging a strict
control $\hat{U}_{t}$, that is $\hat{q}_{t}\left(  du\right)  =\delta_{\hat
{U}_{t}}\left(  du\right)  $.
\end{corollary}

\section{Conclusion}

We have proved two results on the existence of an optimal control for systems
governed by decoupled as well as coupled FBSDEs, by using probabilistic tools.
The ingredients used in the proofs of the main results are based essentially
on tightness techniques on the space $\mathcal{C}$ of continuous functions as
well as on the space $\mathcal{D}$ of c\`{a}dl\`{a}g functions equipped with
Meyer-Zheng topology or Jakubowsky S-topology. The assumptions made on the
coefficients are made to ensure weak convergence of the processes under
consideration and the corresponding cost functionals. However, a serious
difficulty remains in the case where the generator depends on the second
backward variable $Z.$ This difficulty consists in finding a natural
assumption ensuring the tightness of the second backward variable $Z.$ This is
exactly the kind of problems encountered when one deals with weak solutions of
BSDEs and coupled FBSDEs with generators depending on $Z$.

\section{Appendix}

The $S$-topology defined by Jakubowski on the space $\mathbb{D}\left(  \left[
0,T\right]  ;\mathbb{R}^{k}\right)  $ of c\`{a}dl\`{a}g functions is weaker
than the Skorokhod topology and the tightness criteria is the same as for the
Meyer-Zheng topology \cite{MZ}. The topology $S$ arises naturally in limit
theorems for stochastic integrals. Let us give some of its properties:

1) \ If $x_{n}\rightarrow_{S}x_{0}$, then $x_{n}\left(  t\right)
\rightarrow_{S}x_{0}\left(  t\right)  $ for each $t$\ except for a countable set.

2) \ If $x_{n}\left(  t\right)  \rightarrow_{S}x_{0}\left(  t\right)  $ for
each $t$\ in a dense set containing $0$\ and $T$\ and $\left\{  x_{n}\right\}
$ is $S$-relatively compact, then $x_{n}\rightarrow_{S}x_{0}$ (not true for
the convergence in measure).

3) \ We recall (see Meyer and Zheng \cite{MZ} \ and Jakubowski \cite{J1,J2})
that for a familly $\left(  X^{n}\right)  _{n}$ of quasi-martingales on the
probability space $\left(  \Omega,\mathcal{F},\mathcal{F}_{t},P\right)  ,$ the
following condition ensures the tightness of the familly $\left(
X^{n}\right)  _{n}$\ on the space $\mathbb{D}\left(  \left[  0,T\right]  ;%
%TCIMACRO{\U{211d} }%
%BeginExpansion
\mathbb{R}
%EndExpansion
^{k}\right)  $\ endowed with the $S$-topology%
\[
\underset{n}{\sup}\left(  \underset{0\leq t\leq T}{\sup}E\left\vert X_{t}%
^{n}\right\vert +CV\left(  X^{n}\right)  \right)  <\infty,
\]
where, for a quasi-martingale $X$ on $\left(  \Omega,\left\{
\mathcal{\mathcal{F}}_{t}\right\}  _{0\leq t\leq T},P\right)  ,CV\left(
X\right)  $ stands for the conditional variation of $X$\ on $\left[
0,T\right]  $, and is defined by%
\[
CV\left(  X\right)  =\sup E\left(  \underset{i}{\sum}\left\vert E\left(
X_{t_{i}+1}-X_{t_{i}}\right)  \diagup\mathcal{\mathcal{F}}_{t_{i}}%
^{n}\right\vert \right)  ,
\]
where the supremum is taken over all partitions of $\left[  0,T\right]  $.

Let $N^{a,b}\left(  Y\right)  $ denotes the number of up-crossing of the
function $Y\in\mathbb{D}\left(  \left[  0,T\right]  ;\mathbb{R}^{k}\right)  $
in given levels $a<b$ (recall that $N^{a,b}\left(  Y\right)  \geq k$ if one
can find numbers $0\leq t_{1}<t_{2}<\cdot\cdot\cdot<t_{2k-1}<t_{2k}\leq T$
such that $Y\left(  t_{2i-1}\right)  <a$ and $Y\left(  t_{2i}\right)
>b,i=1,2,...,k).$

\begin{lemma}
\label{A4} (A criteria for $S$-tightness). \textit{A sequence }$\left(
Y^{n}\right)  _{n\in\mathbb{N}}$\textit{\ is }$S$\textit{-tight if and only if
it is relatively compact on the }$S$\textit{-topology. Let }$\left(
Y^{n}\right)  _{n\in\mathbb{N}}$\textit{\ be a family of stochastic processes
in }$\mathbb{D}\left(  \left[  0,T\right]  ;\mathbb{R}^{k}\right)  $\textit{.
Then this family is tight for the }$S$ \textit{-topology if and only if
}$\left(  \left\Vert Y^{n}\right\Vert _{\infty}\right)  _{n}$\textit{\ and
}$N^{a,b}\left(  Y^{n}\right)  $\textit{\ are tight for each }$a<b.$
\end{lemma}

\begin{lemma}
\label{A5} (The a.s. Skorokhod representation ).\textit{\ Let }$\left(
\mathbb{D},S\right)  $\textit{\ be a topological space on which there exists a
countable family of }$S$\textit{-continuous functions separating points in
}$X$\textit{. Let }$\left\{  X_{n}\right\}  _{n\in\mathbb{N}}$\textit{\ be a
uniformly tight sequence of laws on }$\mathbb{D}$\textit{. In every
subsequence }$\left\{  X_{n_{k}}\right\}  $\textit{\ one can find a further
subsequence }$\left\{  X_{n_{k_{l}}}\right\}  $\textit{\ and stochastic
processes }$\left\{  Y_{l}\right\}  $\textit{\ defined on }$\left(  \left[
0,T\right]  ,\mathcal{B}_{\left[  0,T\right]  },l\right)  $\textit{\ such
that}
\end{lemma}

\begin{equation}
Y_{l}\sim X_{n_{k_{l}}}, \ l=1,2,... \tag{1}%
\end{equation}

\textit{for each }$w\in\left[  0,T\right]  $
\begin{equation}
Y_{l}\left(  w\right)  \underset{S}{\rightarrow}Y_{0}\left(  w\right)  ,\text{
as }l\rightarrow\infty, \tag{2}%
\end{equation}

\textit{and for each }$\varepsilon>0$\textit{, there exists an }
$S$\textit{-compact subset }$K_{\varepsilon}\subset\mathbb{D}$\textit{\ such
that}%
\begin{equation}
P\left(  \left\{  w\in\left[  0,T\right]  :Y_{l}\left(  w\right)  \in
K_{\varepsilon},l=1,2,...\right\}  \right)  >1-\varepsilon. \tag{3}%
\end{equation}

One can say that (2) and (3) describe "the almost sure convergence in
compacts" and that (1), (2) and (3) define the strong a.s. Skorokhod
representation for subsequences ("strong" because of condition (3)).

\begin{remark}
\label{A6} The projection \ $\pi_{T}:\,\,y\in(\mathbb{D}([0,\,T];\,\mathbb{R}%
),\,S)\longmapsto\,y(T)$ is continuous (see Remark 2.4, p.8 in Jakubowski
\cite{J1}), but $y\longmapsto\,y(t)$ is not continuous for each $0\leq t\leq
T$.
\end{remark}

\begin{lemma}
\label{A7} Let $(X^{n},\,M^{n})$ be a multidimensional process in
$\mathbb{D}([0,\,T];\,\mathbb{R}^{p})\,(p\in\mathbb{R}^{\ast})$ converging to
$(X,\,M)$ in the S-topology. Let $(\mathcal{F}_{t}^{X^{n}})_{t\geq0}$ (resp.
$(\mathcal{F}_{t}^{X})_{t\geq0}$) be the minimal complete admissible
filtration for $X^{n}$ (resp.$X$). We assume that $\sup_{n}E\left[
\sup_{0\leq t\leq T}|M_{t}^{n}|^{2}\right]  <C_{T}\,\,\forall T>0,\,M^{n}$ is
a $\mathcal{F}^{X^{n}}$-martingale and $M$ is a $\mathcal{F}^{X}$-adapted.
Then $M$ is a $\mathcal{F}^{X}$-martingale.
\end{lemma}

\begin{lemma}
\label{A8} Let $(Y^{n})_{n>0}$ be a sequence of processes converging weakly in
$\mathbb{D}([0,\,T];\,\mathbb{R}^{p})$ to $Y$. We assume that $\sup
_{n}E\left[  \sup_{0\leq t\leq T}|Y_{t}^{n}|^{2}\right]  <+\infty$. Then, for
any $t\geq0,\,E\left[  \sup_{0\leq t\leq T}|Y_{t}|^{2}\right]  <+\infty$.
\end{lemma}

\noindent\textbf{Acknowledgments}. 1)A large part of this work has been
carried out when the third author was visiting the Laboratoire LAMAV,
Universit\'{e} de Valenciennes (France) in June 2014. He is grateful for warm
hospitality and support.

2) The authors would like to thank the anonymous referee for very useful
suggestions, which lead to an improvement of the paper.

\end{document}